# A Stability Criterion for Nonparametric Minimal Submanifolds

Yng-Ing Lee and Mu-Tao Wang

October 15, 2002

**Abstract**

An $n$ dimensional minimal submanifold $\Sigma$ of $\mathbb{R}^{n+m}$ is called nonparametric if $\Sigma$ can be represented as the graph of a vector-valued function $f : D \subset \mathbb{R}^n \mapsto \mathbb{R}^m$. This note provides a sufficient condition for the stability of such $\Sigma$ in terms of the norm of the differential $df$.

## 1 Introduction

A minimal submanifold is stable if the second derivative of the volume functional with respect to any compact supported normal variational field is nonnegative. A non-parametric minimal hypersurface in the Euclidean space $\mathbb{R}^{n+1}$ is always stable. This is no longer true when the codimension is greater than one. A non-stable non-parametric minimal surface in four dimension was constructed by Lawson and Osserman in [5]. It seems very little is known about the stability of higher codimension minimal submanifolds except for calibrated ones. Recall a submanifold $\Sigma$ is calibrated by a calibrating form $\Omega$ if $\Omega|_\Sigma$ is the volume form of $\Sigma$, or $*\Omega = 1$ where $*$ is the Hodge star operator. In particular, a non-parametric hypersurface in $\mathbb{R}^{n+1}$ is calibrated by the $n$ form $i(N)dx^1 \cdots dx^{n+1}$ where $N$ is a extension of the unit normal vector field of $\Sigma$.

In [8], the second author constructs solutions to the Dirichlet problem of minimal surface systems in higher dimension and codimension and the solutions satisfies $*\Omega > \frac{1}{2}$ for $\Omega$ the volume form of an $n$-dimensional subspace. When $\Sigma$ is the graph of $f : D \subset \mathbb{R}^n \mapsto \mathbb{R}^m$ and $\Omega$ is the volume form of the domain $\mathbb{R}^n$ extending to the whole $\mathbb{R}^{n+m}$, we have the relation $*\Omega = \frac{1}{\sqrt{\det(I+(df)^T df)}}$. $*\Omega$ is actually the Jacobian of the projection map $\pi : \Sigma \mapsto D$. In particular, a lower bound on $*\Omega$ implies an upper bound on the norm of $df$. In this paper, we discover a criterion for the stability of minimal submanifolds in terms of such condition.

**Theorem A.** *Let $\Sigma$ be the graph of $f : D \subset \mathbb{R}^n \mapsto \mathbb{R}^m$. If $\Sigma$ is minimal and $||df|| \leq \frac{\sqrt{1+c}-1}{\sqrt{c}}$, then $\Sigma$ is stable. Here the norm $||df||$ is defined to be $\sup_{|V|=1} |df(V)|$ and c is a constant which is 1 if $m \leq 2$ or $n \leq 2$, and is $\min(m-1, n-1)$ in other cases.*

Since $*\Omega = \frac{1}{\sqrt{\det(I+(df)^T df)}} = \frac{1}{\sqrt{\prod(1+\lambda_i^2)}}$, it is not hard to see the following consequence.

**Corollary.** *Let c be the constant as in Theorem A. If $\Sigma$ is minimal and $*\Omega \geq \frac{c}{2(c+1-\sqrt{1+c})}$, then $\Sigma$ is stable*

In particular, when $m \leq 2$ or $n \leq 2$, if $\Sigma$ is minimal and $*\Omega \geq \frac{2+\sqrt{2}}{4}$, then $\Sigma$ is stable. The condition $*\Omega > \frac{c}{2(c+1-\sqrt{1+c})}$ corresponds a region in the Grassmannian. For minimal surfaces in $\mathbb{R}^3$, Barbosa and Do Carmo [1],[2] proved if the area of the image of the Gauss map is less than $2\pi$ then $D$ is stable. Their result was also obtained by D. Fischer-Colbrie and R. Schoen in [4] by a different method. This raises the general question of how to characterize stability by the Gauss map of a minimal submanifold. We remark the condition is not sharp in view of the codimension one case.

The proof of Theorem A utilizes a second variation formula of R. Mclean [6] for calibrated submanifolds. The corresponding formula for complex submanifolds of Kähler manifolds was derived by J. Simons [7].

Part of this paper is completed while both authors were visiting the National Center of Theoretical Science in National Tsing-Hua University, Hsin-Chu, Taiwan. The authors wish to express their gratitude for the excellent support provided by the center during their stays. The second author would like to thank Ben Andrews, Bob Gulliver and Brian White for inspring discussions.

## 2  Second Variation Formula

We first recall the second variation formula for minimal submanifolds. Let $F_0 : \Sigma \mapsto \mathbb{R}^{n+m}$ be a minimal submanifold and let $F : \Sigma \times [0,1) \mapsto \mathbb{R}^{n+m}$ be a one-parameter family of immersions with $F(\cdot, 0) = F_0$. We may assume the variation field $V = F_*(\frac{\partial}{\partial s})$ is normal and of compact support. For simplicity, we will identify $F_0(\Sigma)$ with $\Sigma$ and denote $F(\cdot, s)$ by $F_s$. A coordinate system $\{x^i\}$ in a neighborhood of $p \in \Sigma$ is fixed. Let $g_{ij}(s)$ be the induced metric and $dv_s = \sqrt{\det g_{ij}(s)}\, dX$ be the volume form on $F_s(\Sigma)$. At $s=0$, the volume form will be written as $dv$ instead.

We recall the second variation formula from [3]:
$$\frac{d^2}{ds^2}\Big|_{s=0} \int_\Sigma dv_s = \int_\Sigma (||\nabla^N V||^2 - B(V,V))\, dv$$

where $\nabla^N V$ is the covariant derivative of $V$ as a section of the normal bundle and $B(V,V) = \sum_{ijkl} g^{ik} g^{jl} \langle \frac{\partial^2 F}{\partial x^i \partial x^j}, V\rangle \langle \frac{\partial^2 F}{\partial x^k \partial x^l}, V\rangle$. The minimal submanifold $\Sigma$ is stable if and only if

$$\int_\Sigma ||\nabla^N V||^2\, dv \geq \int_\Sigma B(V,V)\, dv$$

for any normal vector field with compact support. Since the second variation formula does not depend on $\frac{\partial^2 F}{\partial s^2}$, we may consider only the case $\frac{\partial^2 F}{\partial s^2}$ is zero at $s=0$. In this case, the following equation holds at every point.

$$\frac{\partial^2}{\partial s^2}\Big|_{s=0} \sqrt{\det g_{ij}(s)} = (||\nabla^N V||^2 - B(V,V))\sqrt{\det g_{ij}} \qquad (2.1)$$

In [6], a different second variation formula is derived in the presence of a calibrating form $\Omega$. In the following, we derive the formula for completeness. We shall assume $\Omega$ is locally an exact form.

Now $\int_{F_s(\Sigma)} \Omega = \int_\Sigma F_s^* \Omega$ is a constant. Write $F_s^* \Omega = *\Omega(s) \sqrt{\det g_{ij}(s)}\, dX$, where $*\Omega = \frac{1}{\sqrt{\det g_{ij}}} \Omega(\frac{\partial F}{\partial x^1}, \cdots \frac{\partial F}{\partial x^n})$.

Since
$$0 = \frac{d^2}{ds^2} \int_\Sigma F_s^* \Omega$$
$$= \int_\Sigma [(\frac{\partial^2}{\partial s^2} *\Omega)\sqrt{\det g_{ij}(s)} + 2(\frac{\partial}{\partial s} *\Omega)(\frac{\partial}{\partial s}\sqrt{\det g_{ij}(s)})$$
$$+ *\Omega \frac{\partial^2}{\partial s^2} \sqrt{\det g_{ij}(s)}]\, dX$$

and $\frac{\partial}{\partial s}|_{s=0}\sqrt{\det g_{ij}}(s) = 0$ by the minimal condition, at $s=0$ we have

$$\int_\Sigma *\Omega \frac{\partial^2}{\partial s^2}\sqrt{\det g_{ij}}(s)\, dX = -\int_\Sigma (\frac{\partial^2}{\partial s^2} *\Omega)\, dv$$

That is,

$$\int_\Sigma *\Omega(||\nabla^N V||^2 - B(V,V))\, dv = -\int_\Sigma (\frac{\partial^2}{\partial s^2}|_{s=0} *\Omega)\, dv \qquad (2.2)$$

We shall compute $\frac{\partial^2}{\partial s^2} *\Omega$ using the formula $*\Omega = \frac{1}{\sqrt{\det g_{ij}}}\Omega(\frac{\partial F}{\partial x^1}, \cdots \frac{\partial F}{\partial x^n})$.
Thus

$$\begin{aligned}
\frac{\partial^2}{\partial s^2} *\Omega &= \frac{\partial^2}{\partial s^2}(\frac{1}{\sqrt{\det g_{ij}}})\Omega(\frac{\partial F}{\partial x^1}, \cdots \frac{\partial F}{\partial x^n}) \\
&+ 2\frac{\partial}{\partial s}(\frac{1}{\sqrt{\det g_{ij}}})\frac{\partial}{\partial s}\Omega(\frac{\partial F}{\partial x^1}, \cdots \frac{\partial F}{\partial x^n}) \\
&+ (\frac{1}{\sqrt{\det g_{ij}}})\frac{\partial^2}{\partial s^2}\Omega(\frac{\partial F}{\partial x^1}, \cdots \frac{\partial F}{\partial x^n})
\end{aligned} \qquad (2.3)$$

At $s=0$, the minimal condition implies the second term vanishes and the first term becomes

$$\frac{\partial^2}{\partial s^2}|_{s=0}(\frac{1}{\sqrt{\det g_{ij}}}) = -(\det g_{ij})^{-1}\frac{\partial^2}{\partial s^2}|_{s=0}(\sqrt{\det g_{ij}})$$

Because $\frac{\partial^2 F}{\partial s^2}$ is zero at $s=0$, the third term is

$$\frac{\partial^2}{\partial s^2}\Omega(\frac{\partial F}{\partial x^1}, \cdots \frac{\partial F}{\partial x^n}) = 2[\Omega(\frac{\partial^2 F}{\partial s \partial x^1}, \frac{\partial^2 F}{\partial s \partial x^2} \cdots \frac{\partial F}{\partial x^n}) + \cdots] \qquad (2.4)$$

Denote $\frac{\partial F}{\partial x^i}$ by $\partial_i$, then $\frac{\partial^2 F}{\partial s \partial x^i} = \frac{\partial^2 F}{\partial x^i \partial s} = \nabla_{\partial_i} V$ where $\nabla$ is the connection on $\mathbb{R}^{n+m}$ and $V = \frac{\partial F}{\partial s}$ is the variation field. In the following computation $(\cdot)^T$ and $(\cdot)^N$ denote the tangent and normal part of vector respectively.

$$\begin{aligned}
&\Omega(\frac{\partial^2 F}{\partial s \partial x^1}, \frac{\partial^2 F}{\partial s \partial x^2} \cdots \frac{\partial F}{\partial x^n}) \\
&= \Omega((\nabla_{\partial_1} V)^T + (\nabla_{\partial_1} V)^N, (\nabla_{\partial_2} V)^T + (\nabla_{\partial_2} V)^N, \cdots \partial_n) \\
&= \Omega((\nabla_{\partial_1} V)^T, (\nabla_{\partial_2} V)^T, \cdots \partial_n) + \Omega((\nabla_{\partial_1} V)^T, (\nabla_{\partial_2} V)^N, \cdots \partial_n) \\
&+ \Omega((\nabla_{\partial_1} V)^N, (\nabla_{\partial_2} V)^T, \cdots \partial_n) + \Omega((\nabla_{\partial_1} V)^N, (\nabla_{\partial_2} V)^N, \cdots \partial_n)
\end{aligned}$$

We can assume $\{x^i\}$ is a normal coordinate system in a neighborhood of $p$ with respect to the induced metric on $\Sigma$. Hence $g_{ij}(0) = \delta_{ij}$ at $p$. We do the computation at point $p$ and get

$$\Omega((\nabla_{\partial_1} V)^T, (\nabla_{\partial_2} V)^T, \cdots \partial_n) = *\Omega(\langle V, \nabla_{\partial_1}\partial_1\rangle\langle V, \nabla_{\partial_2}\partial_2\rangle - \langle V, \nabla_{\partial_1}\partial_2\rangle^2)$$

Continue from equation (2.4), we derive

$$2[\Omega(\frac{\partial^2 F}{\partial s \partial x^1}, \frac{\partial^2 F}{\partial s \partial x^2} \cdots \frac{\partial F}{\partial x^n}) + \cdots]$$
$$= 2 \sum_{i<j} \Omega(\partial_1, \cdots (\nabla_{\partial_i} V)^T, \cdots (\nabla_{\partial_j} V)^N, \cdots \partial_n)$$
$$+ 2 \sum_{i<j} \Omega(\partial_1, \cdots (\nabla_{\partial_i} V)^N, \cdots (\nabla_{\partial_j} V)^T, \cdots \partial_n)$$
$$+ 2 \sum_{i<j} \Omega(\partial_1, \cdots (\nabla_{\partial_i} V)^N, \cdots (\nabla_{\partial_j} V)^N, \cdots \partial_n)$$
$$+ 2 * \Omega(\sum_{i<j} \langle V, \nabla_{\partial_i}\partial_i\rangle\langle V, \nabla_{\partial_j}\partial_j\rangle - \langle V, \nabla_{\partial_i}\partial_j\rangle^2)$$

It thus follows from (2.3) that

$$\frac{\partial^2}{\partial s^2}|_{s=0} * \Omega = - *\Omega ||\nabla^N V||^2 + *\Omega \sum_{i,j} \langle V, \nabla_{\partial_i}\partial_j\rangle^2$$
$$+ 2 \sum_{i<j} \Omega(\partial_1, \cdots (\nabla_{\partial_i} V)^T, \cdots (\nabla_{\partial_j} V)^N, \cdots \partial_n)$$
$$+ 2 \sum_{i<j} \Omega(\partial_1, \cdots (\nabla_{\partial_i} V)^N, \cdots (\nabla_{\partial_j} V)^T, \cdots \partial_n)$$
$$+ 2 \sum_{i<j} \Omega(\partial_1, \cdots (\nabla_{\partial_i} V)^N, \cdots (\nabla_{\partial_j} V)^N, \cdots \partial_n)$$
$$+ 2 * \Omega(\sum_{i<j} \langle V, \nabla_{\partial_i}\partial_i\rangle\langle V, \nabla_{\partial_j}\partial_j\rangle - \langle V, \nabla_{\partial_i}\partial_j\rangle^2)$$

However,

$$*\Omega \sum_{i,j} \langle V, \nabla_{\partial_i}\partial_j\rangle^2 + 2 * \Omega(\sum_{i<j} \langle V, \nabla_{\partial_i}\partial_i\rangle\langle V, \nabla_{\partial_j}\partial_j\rangle - \langle V, \nabla_{\partial_i}\partial_j\rangle^2)$$
$$= *\Omega(\sum_i \langle V, \nabla_{\partial_i}\partial_i\rangle)^2 = 0$$

Therefore,

$$\frac{\partial^2}{\partial s^2}|_{s=0} * \Omega = - *\Omega ||\nabla^N V||^2 + 2\sum_{i<j}\Omega(\partial_1,\cdots(\nabla_{\partial_i}V)^T,\cdots(\nabla_{\partial_j}V)^N,\cdots\partial_n)$$

$$+ 2\sum_{i<j}\Omega(\partial_1,\cdots(\nabla_{\partial_i}V)^N,\cdots(\nabla_{\partial_j}V)^T,\cdots\partial_n)$$

$$+ 2\sum_{i<j}\Omega(\partial_1,\cdots(\nabla_{\partial_i}V)^N,\cdots(\nabla_{\partial_j}V)^N,\cdots\partial_n)$$

Combine equations (2.1) and (2.2), we obtain

**Proposition 2.1** *Let $\Omega$ be an exact parallel n-form and $\Sigma$ be an n dimensional minimal submanifold in $\mathbb{R}^{n+m}$. Assume that $V$ is a normal variation field and $\nabla_V V = 0$ along $\Sigma$. Then one has*

$$\begin{aligned}
\int_\Sigma &*\Omega(||\nabla^N V||^2 - B(V,V))\,dv \\
= \int_\Sigma &[*\Omega||\nabla^N V||^2 \\
&- 2\sum_{i<j}\Omega(\partial_1,\cdots,(\nabla_{\partial_i}V)^T,\cdots,(\nabla_{\partial_j}V)^N,\cdots,\partial_n) \\
&- 2\sum_{i<j}\Omega(\partial_1,\cdots,(\nabla_{\partial_i}V)^N,\cdots,(\nabla_{\partial_j}V)^N,\cdots,\partial_n) \\
&- 2\sum_{i<j}\Omega(\partial_1,\cdots,(\nabla_{\partial_i}V)^T,\cdots,(\nabla_{\partial_j}V)^N,\cdots,\partial_n)]\,dv
\end{aligned} \quad (2.5)$$

# 3 Proof of Theorem A

The idea now is to show the right hand side of equation (2.5) is greater than or equal to

$$\delta \int_\Sigma *\Omega(||\nabla^N V||^2 - B(V,V))\,dv$$

for some $\delta < 1$. We shall express the integrand in the right hand side of equation (2.5) in terms of a particular orthonormal basis. At any point $p$, we consider the singular value decomposition of $df: \mathbb{R}^n \mapsto \mathbb{R}^m$. We have

$$df(a_i) = \lambda_i a_{n+i}$$

where $\lambda_i \geq 0$ are the singular values of $df$, or eigenvalues of $\sqrt{(df)^T df}$. $\{a_i\}_{i=1,\cdots,n}$ is an orthonormal basis of eigenvectors of $\sqrt{(df)^T df}$. The set $\{a_{n+i}\}$ can be completed to form an orthonormal basis $\{a_\alpha\}_{\alpha=n+1,\cdots,n+m}$ for $\mathbb{R}^m$. (In case $m < n$, we will have $\lambda_i = 0$ for $i > m$ and $\{a_{n+i}\}_{i=1,\cdots,m}$ forms an orthonormal basis for $\mathbb{R}^m$.) Now $\{e_i = \frac{1}{\sqrt{1+\lambda_i^2}}(a_i + \lambda_i a_{n+i})\}$ and $\{e_{n+i} = \frac{1}{\sqrt{1+\lambda_i^2}}(a_{n+i} - \lambda_i a_i)\}$ can be completed to give orthonormal basis of the tangent and normal space. The orthonormal basis for the normal space is denoted by $\{e_\alpha\}_{\alpha=n+1,\cdots,n+m}$. In these bases we denote $(\nabla_{e_i} V)^N = \sum_\alpha V_i^\alpha e_\alpha$ and $(\nabla_{e_i} V)^T = -\sum_{\alpha,j} V^\alpha h_{\alpha ij} e_j$. We shall assume $m \geq n$ in the following calculation, the other case can be treated similarly.

At the point $p$ in these bases, the integrand of the right hand side of equation (2.5) can be written as

$$*\Omega[\sum_{i,\alpha}(V_i^\alpha)^2 - 2\sum_{i<j}\lambda_i\lambda_j V_i^{n+i}V_j^{n+j} + 2\sum_{i<j}\lambda_i\lambda_j V_i^{n+j}V_j^{n+i}$$
$$+ 2\sum_\alpha\sum_{i<j} V^\alpha h_{\alpha ii} V_j^{n+j}\lambda_j - 2\sum_\alpha\sum_{i<j} V^\alpha h_{\alpha ij} V_j^{n+i}\lambda_i$$
$$+ 2\sum_\alpha\sum_{i<j} V^\alpha h_{\alpha jj} V_i^{n+i}\lambda_i - 2\sum_\alpha\sum_{i<j} V^\alpha h_{\alpha ij} V_i^{n+j}\lambda_j]$$

which is the same as

$$*\Omega[\sum_{i,\alpha}(V_i^\alpha)^2 - \sum_{i\neq j}\lambda_i\lambda_j V_i^{n+i}V_j^{n+j} + \sum_{i\neq j}\lambda_i\lambda_j V_i^{n+j}V_j^{n+i}$$
$$+ 2\sum_\alpha\sum_{i\neq j} V^\alpha h_{\alpha ii} V_j^{n+j}\lambda_j - 2\sum_\alpha\sum_{i\neq j} V^\alpha h_{\alpha ij} V_j^{n+i}\lambda_i]$$

By minimality, we have

$$2\sum_{i\neq j} V^\alpha h_{\alpha ii} V_j^{n+j}\lambda_j = -2\sum_j V^\alpha h_{\alpha jj} V_j^{n+j}\lambda_j$$

Define $\Xi$ by

$$\Xi = \sum_{i,\alpha}(V_i^\alpha)^2 - \sum_{i\neq j}\lambda_i\lambda_j V_i^{n+i}V_j^{n+j} + \sum_{i\neq j}\lambda_i\lambda_j V_i^{n+j}V_j^{n+i}$$
$$- 2\sum_\alpha\sum_{i,j} V^\alpha h_{\alpha ij} V_j^{n+i}\lambda_i$$

If we can show

$$\Xi \geq \delta[\sum_{i,\alpha}(V_i^\alpha)^2 - \sum_{ij}(\sum_\alpha V^\alpha h_{\alpha ij})^2]$$

for some $\delta < 1$, then we are done in review of equation (2.5).

By Cauchy-Schwarz inequality, for any $\epsilon > 0$ to be determined,

$$-2\sum_{i,j,\alpha} V^\alpha h_{\alpha ij} V_j^{n+i}\lambda_i \geq -\epsilon \sum_{i,j}(\sum_\alpha V^\alpha h_{\alpha ij})^2 - \frac{1}{\epsilon}\sum_{i,j}(V_j^{n+i}\lambda_i)^2$$

Therefore

$$\Xi \geq \sum_{i\alpha}(V_i^\alpha)^2 - \sum_{i\neq j}\lambda_i\lambda_j V_i^{n+i}V_j^{n+j} + \sum_{i\neq j}\lambda_i\lambda_j V_i^{n+j}V_j^{n+i} - \frac{1}{\epsilon}\sum_{i,j}(V_j^{n+i}\lambda_i)^2$$
$$- \epsilon \sum_{i,j}(\sum_\alpha V^\alpha h_{\alpha ij})^2$$

Now assume each $\lambda_i \leq \eta$, then

$$\Xi \geq \sum_{i,\alpha}(V_i^\alpha)^2 - \eta^2\sum_{i\neq j}|V_i^{n+i}||V_j^{n+j}| - \eta^2\sum_{i\neq j}|V_i^{n+j}||V_j^{n+i}| - \frac{\eta^2}{\epsilon}\sum_{i,j}(V_j^{n+i})^2$$
$$- \epsilon\sum_{i,j}(\sum_\alpha V^\alpha h_{\alpha ij})^2 \quad (3.1)$$

By Cauchy-Schwarz inequality

$$\sum_{i\neq j}|V_i^{n+i}||V_j^{n+j}| \leq (n-1)\sum_i (V_i^{n+i})^2$$

and

$$\sum_{i\neq j}|V_i^{n+j}||V_j^{n+i}| \leq \sum_{i\neq j}(V_i^{n+j})^2$$

In a general case, the coefficient in the right hand side of the first inequality should be $\min(m-1, n-1)$. Let $c = 1$ if $m \leq 2$ or $n \leq 2$ and $c = \min(m-1, n-1)$ in other cases. Plug into equation (3.1), we obtain,

$$\Xi \geq (1 - \frac{\eta^2}{\epsilon} - c\eta^2)\sum_{i,\alpha}(V_i^\alpha)^2 - \epsilon\sum_{i,j}(\sum_\alpha V^\alpha h_{\alpha ij})^2 \quad (3.2)$$

We need to have $1 - \frac{\eta^2}{\epsilon} - c\eta^2 \geq \epsilon$ which is $\eta^2 \leq \frac{\epsilon(1-\epsilon)}{1+c\epsilon}$. It is not hard to see $\max_{0<\epsilon<1} \frac{\epsilon(1-\epsilon)}{1+c\epsilon} = \frac{(\sqrt{c+1}-1)^2}{c}$ when $\epsilon = \frac{\sqrt{1+c}-1}{c}$. Thus if we assume $||df|| \leq \frac{\sqrt{1+c}-1}{\sqrt{c}}$, then each $|\lambda_i| \leq \frac{\sqrt{1+c}-1}{\sqrt{c}}$ and

$$\Xi \geq \frac{\sqrt{1+c}-1}{c}[\sum_{i,\alpha}(V_i^\alpha)^2 - \sum_{i,j}(\sum_\alpha V^\alpha h_{\alpha ij})^2]$$

Theorem A is proved.

## 4 Examples

The construction in [8] supplies examples for such stable minimal submanifolds. Given any $\phi : D \subset \mathbb{R}^n \mapsto \mathbb{R}^m$ defined on a convex domain $D$, we can scale $\phi$ so that on the graph of $\phi$, $*\Omega \geq \frac{c}{2(c+1-\sqrt{1+c})}$ and the derivative of $\phi$ satisfies the requirement in [8]. It was proved in [8] that the Cauchy-Dirichlet problem of the mean curvature flow for initial data $\phi$ is solvable and $*\Omega \geq \frac{c}{2(c+1-\sqrt{1+c})}$ is preserved along the flow. The flows converges to a minimal submanifold which is stable by the Corollary.

## References

‑